\documentclass[12pt,leqno]{amsart}
\usepackage{amssymb,verbatim,enumerate,ifthen,cite}
\usepackage[mathscr]{eucal}
\usepackage{cite}
\usepackage[utf8]{inputenc}
\usepackage[T1]{fontenc}
\def\N{\mathbb{N}}
\def\R{\mathbb{R}}

\def\H{\mathscr{H}}

\newtheorem{theorem}{Theorem}[section]
\newtheorem*{theorem*}{Theorem}
\def\Thm#1#2{\ifthenelse{\equal{#1}{*}}{\begin{theorem*}#2\end{theorem*}}
	{\begin{theorem}\label{T#1}#2\end{theorem}}}

\newtheorem{proposition}[theorem]{Proposition}
\newtheorem*{proposition*}{Proposition}
\def\Prp#1#2{\ifthenelse{\equal{#1}{*}}{\begin{proposition*}#2\end{proposition*}}
	{\begin{proposition}\label{P#1}#2\end{proposition}}}

\newtheorem{corollary}[theorem]{Corollary}
\newtheorem*{corollary*}{Corollary}
\def\Cor#1#2{\ifthenelse{\equal{#1}{*}}{\begin{corollary*}#2\end{corollary*}}
	{\begin{corollary}\label{C#1}#2\end{corollary}}}

\newtheorem{lemma}[theorem]{Lemma}
\newtheorem*{lemma*}{Lemma}
\def\Lem#1#2{\ifthenelse{\equal{#1}{*}}{\begin{lemma*}#2\end{lemma*}}
	{\begin{lemma}\label{L#1}#2\end{lemma}}}

\theoremstyle{definition}
\newtheorem{remark}[theorem]{Remark}
\newtheorem*{remark*}{Remark}
\def\Rem#1#2{\ifthenelse{\equal{#1}{*}}{\begin{remark}\rm #2\end{remark}}
	{\begin{remark}\label{R#1}\rm #2\end{remark}}}

\newtheorem{example}[theorem]{Example}
\newtheorem*{example*}{Example}
\def\Exa#1#2{\ifthenelse{\equal{#1}{*}}{\begin{example*}\rm #2\end{example*}}
	{\begin{example}\label{Ex#1}\rm #2\end{example}}}

\def\eq#1{{\rm(\ref{E#1})}}
\def\Eq#1#2{\ifthenelse{\equal{#1}{*}}
	{\begin{equation*}\begin{aligned}[]#2\end{aligned}\end{equation*}}
	{\begin{equation}\begin{aligned}[]\label{E#1}#2\end{aligned}\end{equation}}}

\begin{document}
	\begin{flushright}
	\end{flushright}
	\vspace{5mm}
	
	\date{\today}
	
\title[Some Stability Results on Graphs]
{Some Stability Results on Graphs}

\author[A. R. Goswami]{Angshuman R. Goswami}
\address[A. R. Goswami]{Department of Mathematics, Faculty of Technical Informatics, University of Pannonia,
	H-8200 Veszprem, Hungary}
\email{goswami.angshuman.robin@mik.uni-pannon.hu}
\author[M. K. Shihab]{Mahmood K. Shihab}
\address[M. K. Shihab]{Department of Mathematics, College of Education for Pure Sciences, University of Kirkuk, 36001 Kirkuk, Iraq}
\email{mahmoodkamil30@uokirkuk.edu.iq}


\thanks{The research of the first author was supported by the `Ipar a Veszprémi Mérnökképzésért' foundation.
}

\begin{abstract}
The primary objective of this paper is to introduce Hyers-Ulam-type stability results for monotone, subadditive, and convex graphs. We consider their standard definitions in an approximate sense and demonstrate the existence of a corresponding graph with the same vertex and edge sets bearing the exact ideal structural property. We prove that the weight difference on the two graphs depends on the associated error and does not vary significantly.	
\end{abstract}

\maketitle

\section*{Introduction}
In classical real function theory, stability problems study whether a function that satisfies a given property approximately is close enough in an appropriate sense to a function that satisfies the property exactly, and they seek explicit bounds for the distance between the two functions.\\

The origin of this theory goes back to a question posed by S. M. Ulam in 1940, concerning the stability of homomorphisms between groups. In response to this problem, D. H. Hyers provided the first affirmative answer in 1941. In his seminal paper, Hyers proved that if a mapping between Banach spaces approximately satisfies the additive Cauchy functional equation, then there exists an exact additive mapping near it. This result is now known as the Hyers–Ulam stability theorem.\\

As a simple consequence of Hyers’ result, if a function
$f:\mathbb{R}\rightarrow \mathbb{R}$
satisfies the functional inequality of the form 
\begin{equation*}
	|f(x+y)-f(x)-f(y)|\leq \varepsilon \qquad \qquad (x,y\in\R),
\end{equation*}
then there exists an additive function $g:I\rightarrow \mathbb{R}$ such that  $\|f-g\|\leq \varepsilon/2$ holds (see\cite{Ulam}). This illustrates the fundamental principle of Hyers-Ulam stability: an approximately additive function is uniformly close to an exact additive function.\\

Similar versions of this result were later extended for monotone, convex, subadditive, and many other important classes of functions. More insights in this direction can be found in the books and articles \cite{Hyers, Them, Pales, Goswami, Hardy}. In the last few decades, mathematicians studied the idea of Hyers-Ulam-type stability for functions involving ordinary/ partial differential equations and inequalities. Some of the well-recognized works are included in the papers \cite{Obloza,Jung, Miura,Takahasi,Brzdek} and references therein. Motivated by all this, we conducted some stability-oriented investigations in sequences. Representative results on this subject are included in \cite{Goswamii, Goswamiii}. Now, we aim to extend our studies by inspecting and integrating the possible underlying Hyers-Ulam-type stability characteristics in a graph setting.\\

In this paper, we consider three obvious properties associated with classical graphs, namely monotonicity, subadditivity, and convexity. Various parameters of a graph such as maximum degree($\Delta(G)$), chromatic number($\chi(G)$), clique number($w(G)$), independence number($\alpha(G)$), number of connected components, treewidth, etc. bear these structural characteristics and naturally leads to the definitions of monotone, subadditive, and convex graphs. A number of well-known and influential results concerning these parameters are available in \cite{Jakovac,Saleh, Dirac, Rezaei, Zhao, Schaefer, Brooks, Erdos, Turan, Vizing,Renyi,Devroye}. Before we discuss the terminologies and present our findings, we need to recall some basic notations.\\

Throughout the paper, $\N$ and $\R_+$ denote the set of natural and positive real numbers, respectively. The symbol $G(V,E)$ is used to denote the non-weighted graph $G$ along with its vertex and edge sets $V$ and $E$ respectively. We assume $\H$ is the collection of all non-empty subsets of $G$ and based on it we define a mapping $w:\H\to\R_+$. In other words, we assign weights to each and every subgraph of $G$, and hence, we can term $G(V,E,w)$ as a weighted graph. It is also important to note that despite of having the same vertex and edge sets, the two graphs $G(V,E,w_1)$ and $G(V,E,w_2)$ are not identical as they bear two distinct weights. Another important observation is that $(\H,\subseteq)$ is a poset, and $G$ is its greatest element.\\

For a fixed $\varepsilon>0$, a weighted graph $G(V,E,w)$ is called \textit{approximately monotone}, if for any \newline
$H,H'\in\H$ satisfy the following discrete functional inequality

\Eq{1}{
	w(H)\leq w(H')+\varepsilon \qquad\mbox{provided}\qquad H\subset H'\,.
}
Clearly, if $\varepsilon=0$, then $G(V,E,w)$ is termed a \textit{monotone(increasing)} graph. We establish that if any weighted graph $G$ satisfies the inequality \eq{1}, then there exists an ordinary monotone graph $G(V,E,\widetilde{w})$, such that $\|w-\widetilde{w}\|\leq\varepsilon˛/2$ holds for all $H\in\H$. The converse assertion is also validated.\\

Let $\varepsilon>0$ be fixed. A weighted graph $G(V,E,w)$ is said to be \textit{approximately subadditive}, if for any  
$H_1,\cdots ,H_n\in\H$, the weight function $w$ satisfies the following discrete inequality

\Eq{2}{
w(H)\leq w(H_1)+\cdots+w(H_n)+\varepsilon \qquad\mbox{where}\qquad H=\overset{n}{\underset{i=1}\cup} H_i\,.
}
It is evident that for $\varepsilon=0$, the weighted graph  $G(V,E,w)$ is nothing but a \textit{subadditive graph}. We demonstrate that if a graph $G(V,E,w)$ is approximately subadditive, then there exists a usual subadditive minorant $G(V,E,\overline{w})$ satisfying the condition $w-\overline{w}\leq\varepsilon$.\\

Since convexity is also referred to as second-order monotonicity, we construct a meaningful definition of a convex graph by taking inspiration from the theory of sequential convexity (see \cite{Mitrinovicc}). A weighted graph $G(V,E,{w})$ is said to be \textit{convex} if any arbitrarily chosen three subgraphs  $\underline{H}$,$H$,$\overline{H}\in\H$ with $\underline{H}\subseteq H\subseteq \overline{H}$, are equipped with the following discrete functional inequality
\Eq{*}{
	2w(H)\leq w(\underline{H})+w(\overline{H})\,.
}
The approximately convex graph can be defined in an analogous way, for which we also propose a decomposition result showing Hyers-Ulam-type stability.\\

We start our investigation with approximately monotone graphs.

\section{Main Results}
To maintain full generality, we have not imposed any restrictions on the cardinalities of the vertex and edge sets $V$ and $E$ of $G$. Nevertheless, the results can be readily adapted to finite graphs, leading naturally to analogous formulations of Hyers-Ulam-type stability.
\Prp{p}{
Let $G(V,E,w)$ be an approximately monotone graph, then there exists a monotone(increasing) graph $G(V,E,\widetilde{w})$ such that $\|w-\widetilde{w}\|\leq \varepsilon/2$ holds. Conversely if $G(V,E,\widetilde{w})$ is a monotone graph satisfying the inequality $\|w-\widetilde{w}\|\leq \varepsilon/2$, then $G(V,E,w)$ is approximately monotone.
}
\begin{proof}
Since $G(V,E,w)$ is approximately monotone, it satisfies the inequality \eq{1}, which can also be re-written as follows
\Eq{11}{
w(H)-\dfrac{\varepsilon}{2}\leq w(H')+\dfrac{\varepsilon}{2} \qquad\mbox{provided}\qquad H\subseteq H'\qquad (H,H'\in\H)\,\,.
}
Now we define the function $\widetilde{w}:\H\to\R_+$ as follows
\Eq{*}{
\widetilde{w}(H):=\inf \big\{w(H')\,\,\big|\,\, H\subseteq H'\big\}+\frac{\varepsilon}{2}.
}
Using this definition in \eq{11}, we have

\Eq{*}{
w(H)-\dfrac{\varepsilon}{2}\leq \widetilde{w}(H)\leq w(H)+\dfrac{\varepsilon}{2}\quad\mbox{for all}\quad H\in\H.
}

This shows the norm inequality $\|w-\widetilde{w}\|\leq \dfrac{\varepsilon}{2}$. \\

To establish the monotonicity of $G(V,E,\widetilde{w})$, we assume $H_1,H_2\in\H$ with $H_1\subset H_2$. By using definition $\widetilde{w}$, we can proceed as follows
\Eq{*}{
\widetilde{w}(H_1)=\inf \big\{w(H')\,\,\big|\,\, H_1\subseteq H'\big\}+\frac{\varepsilon}{2}\leq  \inf \big\{w(H')\,\,\big|\,\, H_2\subseteq H'\big\}+\frac{\varepsilon}{2}=\widetilde{w}(H_2).
}
This validates that $G(V,E,\widetilde{w})$ possesses monotonicity and completes the proof of first assertion.\\

To establish the converse assertion, we arbitrarily choose $H,H'\in\H$ such that $H\subset H'$ holds. Now starting with the norm inequality and then using the monotonicity of $G(V,E,\widetilde{w})$, we can compute the following inequality
\Eq{*}{
w(H)\leq \widetilde{w}(H)+\dfrac{\varepsilon}{2}
\leq \widetilde{w}(H')+\dfrac{\varepsilon}{2}
\leq \bigg({w}(H')+\dfrac{\varepsilon}{2}\bigg)+\dfrac{\varepsilon}{2}\leq {w}(H')+{\varepsilon}. 
}
Since $H$ and $H'$ are arbitrarily chosen, the above inequalities represent the inequality \eq{1} showing the approximate monotonicity of $G(V,E,{w})$. This completes the proof.\\
\end{proof}

Now, we discuss subadditive graphs in an approximate sense. Notably, a decreasing graph $G(V,E,w)$ satisfies the property of subadditivity. To verify this claim, suppose that
$G(V,E,w)$ is decreasing and let $H_1,H_2\in\H$ be arbitrary. Then we have
\Eq{*}{
w(H_1\cup H_2)\leq w(H_1)\qquad\mbox{and}\qquad w(H_1\cup H_2)\leq w(H_2).
}
From these two inequalities, it follows immediately that
\Eq{*}{
w(H_1\cup H_2)\leq w(H_1)+w(H_2),
} 
which proves that besides decreasingness $G(V,E,w)$ possesses subadditivity. 
\Thm{q}{
Let $G(V,E,w)$ be an approximately subadditive graph. Then there exists a subadditive graph $G(V,E,\overline{w})$ such that the norm inequality
$\|w-\overline{w}\|<{\varepsilon}$ holds.}
\begin{proof}
For readability, we divide the proof into three parts, each portion establishing one assertion. In the first part, we construct a subadditive minorant $\overline{w}:\H\to\R_+$ of $w$. Next, we prove a sandwich-type result for weighted graphs. Finally, in the third step, combining these two ingredients, we complete the proof of the theorem. \\

\textbf{First Step : }
"Let $G(V,E,w)$ be a weighted graph. We define $\overline{w}:\H\to\R_+$ as follows
\Eq{33}{
\overline{w}\big(H\big):=\inf\Big\{w\big(H_1\big)+\cdots +w\big(H_n\big)\,\,\Big|\,\,H_1,\cdots,H_n\subseteq H\,\,\mbox{such that}\,\,\overset{n}{\underset{i=1}{\cup}}H_i=H \}.
}
Then $G(V,E,\overline{w})$ is a subadditive minorant of $G(V,E,{w})$."\\

From definition, $\overline{w}\big(H\big)\leq w(H)$ is evident. This implies $G(V,E,\overline{w})$ is a minorant of $G(V,E,{w})$.\\

To prove the subadditivity of $G(V,E,\overline{w})$, we assume $H_1, H_2\in \H$ are two distinct subgraphs of G. Then for any $\delta>0$, there exists a partition $H_1^{1},\cdots,H_1^{m}$ for $H_1$ and a another partition $H_2^{1},\cdots,H_2^{n}$ for $H_2$ such that the following two inequalities are satisfied
\Eq{44}{
w\big(H_1^{1}\big)+\cdots +w\big(&H_1^{m}\big)<\overline{w}(H_1)+\dfrac{\delta}{2}\\
&\mbox{and}\\
w\big(H_2^{1}\big)+\cdots +w\big(&H_2^{n}\big)<\overline{w}(H_2)+\dfrac{\delta}{2}.
}
Naturally, the combined partitions of $H_1$ and $H_2$ also act as a partition of $H_1\cup H_2$. Mathematically, this can be represented as follows
\Eq{*}{H_1^1\cup\cdots\cup H_1^m\cup H_2^1\cup\cdots\cup H_2^n=H_1\cup H_2\subseteq G.
}
From \eq{33} and \eq{44}, we obtain the following
\Eq{*}{
\overline{w}\big(H_1\cup H_2\big)&\leq w\big(H_1^{1}\big)+\cdots +w\big(H_1^{m}\big)+w\big(H_2^{1}\big)+\cdots +w\big(H_2^{n}\big)\\
&<\overline{w}(H_1)+\dfrac{\delta}{2}+\overline{w}(H_2)+\dfrac{\delta}{2}\\
&=\overline{w}(H_1)+\overline{w}(H_2)+\delta.
}
Taking $\delta\to 0$, the above inequality can be expressed as
\Eq{*}{
\overline{w}\big(H_1\cup H_2\big)\leq \overline{w}(H_1)+\overline{w}(H_2).
}
Since $H_1,H_2\in\H$ are arbitrary, this yields that the weighted graph $G(V,E,\overline{w})$ possesses subadditivity. This proves the first assertion\\

\textbf{ Second Step : } ``Let $G(V,E,w_1)$ and $G(V,E,w_2)$ be two weighted graphs such that 
\Eq{55}{w_1\big(H\big)\leq w_2\big(H_1\big)+\cdots+w_2\big(H_n\big)
}
holds for any $H_1,\cdots,H_n\subseteq H$ with $\overset{n}{\underset{i=1}{\cup}}H_i=H$.
Then there exists a subadditive graph $G(V,E,\overline{w})$, that satisfies the following discrete functional inequality
\Eq{66}{
w_1\big(H\big)\leq \overline{w}\big(H\big)\leq w_2\big(H\big)\qquad\mbox{for all}\qquad H\in\H.\,\,"}

Just by replacing $w$ with $w_2$ in the right side of \eq{33}, we formulate the required subadditive weight function $\overline{w}:\H\to\R_+$. Moreover, from the inequality \eq{55}, we can observe the following
\Eq{*}{
w_1(H)\leq \inf\Big\{w_2\big(H_1\big)+\cdots +w_2\big(H_n\big)\,\,\Big|\,\,H_1,\cdots,H_n\subseteq H\,\,\mbox{such that}\,\,\overset{n}{\underset{i=1}{\cup}}H_i=H\Big \}\leq w_2(H).
}
This establishes the inequality \eq{66} and validates the complete statement.\\

\textbf{ Third Step : } Finally, we have all the required tools to prove the proposed theorem. Since, $G(V,E,w)$ is approximately subadditive, it satisfies the inequality \eq{2}, which can also be represented in the following way
\Eq{77}{
	w(H)-{\varepsilon}\leq w(H_1)+\cdots+w(H_n) \qquad\mbox{where}\qquad H=\overset{n}{\underset{i=1}{\cup}}H_i\,\,.
}
To deduce a refinement of the above inequality, we define the weight function $w_1:\H\to\R_+$ as 
\Eq{*}{
w_1(H):=\max\big\{0,w(H)-\varepsilon\big\}.
}
Applying $w_1$, we obtain an improved alternative inequality of \eq{77} which is given below
\Eq{88}{
w(H)-{\varepsilon}\leq w_1(H)\leq w(H_1)+\cdots+w(H_n) \qquad\mbox{where}\qquad H=\overset{n}{\underset{i=1}{\cup}}H_i\,\,.
}
As shown earlier, the right-most inequality of \eq{88} ensures the existence of a subadditive weight function $\overline{w}:\H\to\R_+$ sandwiched in between $w_1$ and $w$. Using this result in the inequality \eq{88}, we conclude
\Eq{*}{
w(H)-{\varepsilon}\leq \overline{w}(H)\leq w(H) \qquad \mbox{for all}\quad H\in\H.
}
This shows $\|w-\overline{w}\|\leq \varepsilon$ and completes the proof.\\
\end{proof}

We have already discussed the convexity of graphs in the Introduction section. One of the salient features of a convex graph $G(V,E,w)$ is that if we consider any chain of $(\H,\subseteq)$, then the existence of a global minimum with respect to the corresponding weights is guaranteed. To prove the assertion, we consider a monotone(increasing) chain $C:=\big\{H_n:n\in\N\cup\{0\}\big\}$ of $(\H,\subseteq)$. If the corresponding sequence $\Big(w(H_n)\Big)$ ($n\in\N\cup\{0\})$ also possesses monotonicity, then the assertion is obvious. If the sequence $\Big(w(H_n)\Big)$ is non-monotone, it is sufficient to show that there exists an $H_i\in C$ such that the corresponding weight distributions on the chains
\Eq{*}{
C_{\underline{H_i}}:=\big\{H_0,H_1,\cdots,H_i\big\} \quad &\mbox{is monotonically decreasing}
\\
&\mbox{and}
\\
C_{\overline{H_i}}:={H_i}\cup\Big(C\setminus C_{\underline{H_i}}\Big) \quad &\mbox{is monotonically increasing.}
}
The weighted convexity on $G$ also implies that the sequence $\Big(w(H_n)-w(H_{n-1})\Big)_{n\in\N}$ is monotonically increasing. Therefore there exists a $H_i\in\H$ such that the sub-sequence $\Big(w(H_n)-w(H_{n-1})\Big)_{n=1}^{i}$ is non-positive, while the other sub-sequence $\Big(w(H_n)-w(H_{n-1})\Big)_{n=i}^{\infty}$
possesses non-negativity. Combining these two observations together with their increasingness validates our claim. Next, we propose a stability result for a convex graph, but first, we need to define approximate convexity on graphs.
\\

For a fixed $\varepsilon>0$, a weighted graph $G(V,E,w)$ is called \textit{approximately convex}, if for any $\underline{H}, H, \overline{H}\in\H$ with 
$\underline{H}\subseteq H \subseteq \overline{H}$, the following discrete functional inequality holds
\Eq{99}{
2w(H)\leq w(\underline{H})+w(\overline{H})+\varepsilon\,\,.
}
\Thm{r}
{Let graph $G(V,E,w)$ be an approximately convex graph. Then there exists a convex graph $G(V,E,\widehat{w})$ such that following inequality satisfies
\Eq{369}{
\inf_{H\in\H}w(H)-\dfrac{\varepsilon}{2}\leq \widehat{w}(H)\leq w(H)+\dfrac{\varepsilon}{2}.
}
}
\begin{proof}
To prove this theorem, we first need to establish the following statement:\\
"For any graph $G(V,E,w)$, there exists a non-trivial convex graph $G(V,E,\widehat{w})$ such that the inequality $w(H)\leq \widehat{w}(H)$ holds for all $H\in\H.$"\\

To validate the assertion, for an arbitrary $H\in\H$ we construct the non-negative sequence $\Big(w_k(H)\Big)_{k=0}^{\infty}$ as follows
\Eq{119}{
w_0(H):=\inf_{\underline{H}\subseteq H\subseteq\overline{H}} \dfrac{w(\underline{H})+w(\overline{H})}{2}\qquad\mbox{and}\qquad  w_k(H):=\inf_{\underline{H}\subseteq H\subseteq\overline{H}} \dfrac{w_{k-1}(\underline{H})+w_{k-1}(\overline{H})}{2}\quad \quad (k\in\N).
}
Since the sequence $\Big(w_k(H)\Big)_{k=0}^{\infty}$ is decreasing and bounded below by $0$, this guarantees the convergence of the sequence. Now, we define the weight function $\widehat{w}:\H\to\R_+$ as
\Eq{*}{
  \widehat{w}(H)=\lim_{k\to\infty} w_k(H).
}
Clearly $\widehat{w}\leq w$. Therefore, it will be sufficient to show the convexity of the weight function $\widehat{w}$.\\

Let $H\in\H$ be fixed, and $H_1, H_2\in\H$ are arbitrarily chosen such that $H_1\subseteq H\subseteq H_2$ holds. Then by definition of $\widehat{w}$, we can proceed as follows
\Eq{*}{
\widehat{w}(H)
&= \lim_{k\to\infty} w_k(H)\\
&= \lim_{k\to\infty}\Bigg(\inf_{\underline{H}\subseteq H\subseteq\overline{H}} \dfrac{w_{k-1}(\underline{H})+w_{k-1}(\overline{H})}{2}\Bigg)\\
&\leq\lim_{k\to\infty}\Bigg(\frac{w_{k-1}(H_1)+w_{k-1}(H_2)}{2}\Bigg)\\
&\leq \dfrac{\widehat{w}(H_1)+\widehat{w}(H_2)}{2}.
}
This proves that $G(V,E,\widehat{w})$ is a non-trivial convex minorant of the graph $G(V,E,{w})$.
\\

Now, we have the required ingredient to establish the theorem. We start by re-structuring the inequality \eq{99} as below
\Eq{*}{
\max\bigg\{w(H)-\dfrac{\varepsilon}{2}\,,\,0\bigg\}
\leq \inf_{\underline{H}\subseteq H\subseteq\overline{H}}\dfrac{\Big(w(\underline{H})+\varepsilon/2\Big)+\Big(w(\overline{H})+\varepsilon/2\Big)}{2}\,\,.
}
For simplicity, we use substitutions to represent the above inequality as follows
\Eq{130}{
w'(H)\leq w''(H)\qquad\mbox{for all}\qquad H\in\H.
}\\
where the weight functions $w', w'':\H\to\R_+$ are given by 
\Eq{212}{
w'(H):=\max\bigg\{w(H)-\dfrac{\varepsilon}{2}\,,\,0\bigg\}\qquad\mbox{and}\qquad w''(H):=\inf_{\underline{H}\subseteq H\subseteq\overline{H}}\dfrac{\Big(w(\underline{H})+\varepsilon/2\Big)+\Big(w(\overline{H})+\varepsilon/2\Big)}{2}.
}
Now, we define the sequence $\Big(w''_k\Big)_{k=0}^{\infty}$ as follows
\Eq{1119}{
w''_0(H):=w''(H)\qquad\mbox{and}\qquad  w''_k(H):=\inf_{\underline{H}\subseteq H\subseteq\overline{H}} \dfrac{w''_{k-1}(\underline{H})+w''_{k-1}(\overline{H})}{2}\quad \quad (k\in\N).
}
The construction of the sequence $\Big(w''_k\Big)_{k=0}^{\infty}$ is analogous to \eq{119} and hence, following similar mathematical arguments, we can show the existence of a non-trivial convex weight function $\widehat{w}:\H\to\R_+$ such that for all 
$H\in\H$, the following inequality holds
\Eq{234}{
\qquad\widehat{w}(H)\leq w''_{k}(H)\leq w''(H)\qquad (k\in\N).
}
Also from the inequality \eq{130} and \eq{1119} for all $H\in\H$, we can conclude the following inequality
\Eq{120}{
\inf_{H\in\H}w'(H)\leq \widehat{w}(H)\leq w''_k(H)\qquad (k\in\N).
}
If not, then for a $H\in\H$ there exists a $k\in\N$ such that $w''_k(H)\leq\underset{{H\in\H}}{\inf}w'(H)$ holds. In other words, this implies the validity of the following inequality
\Eq{*}{
w''_k(H)=\inf_{\underline{H}\subseteq H\subseteq\overline{H}} \dfrac{w''_{k-1}(\underline{H})+w''_{k-1}(\overline{H})}{2}\leq\inf_{H\in\H}w'(H).
}
Which yields the existence of either an $\underline{H}\subseteq H$ or an $\overline{H}\supseteq H$ such that at least one of the following inequality satisfied
\Eq{*}{
w''_{k-1}\big(\underline{H}\big)\leq\inf_{H\in\H}w'(H)\qquad\mbox{or}\qquad
w''_{k-1}\big(\overline{H}\big)\leq\inf_{H\in\H}w'(H).
}
Without loss of generality, we assume that from these two inequalities, the first inequality, that is $w''_{k-1}\big(\underline{H}\big)\leq\underset{{H\in\H}}\inf w'(H)$ holds along with the assumed conditions. Now using the same logic, we proceed, and eventually we get the existence of an $H_0\in\H$ that satisfies the following inequality 
\Eq{*}{
w(H_0)\leq\inf_{H\in\H}w'(H).
}
Clearly, this is a contradiction reflecting the inaccuracy of our assumption. Hence, \eq{120} is valid.\\

By combining the two inequalities \eq{234} and \eq{120}, we arrive at 
\Eq{700}{
\inf_{H\in\H}w'(H)\leq \widehat{w}(H)\leq w''(H).
}
Also from \eq{212}, we obtain the following two extensions
\Eq{*}{
\inf_{H\in\H}w(H)-\dfrac{\varepsilon}{2}
&\leq \inf_{H\in\H}\Bigg(\max\bigg\{ w(H)-\dfrac{\varepsilon}{2}\,,0\,\bigg\}\Bigg)=\inf_{H\in\H}w'(H)\\
&\qquad\qquad\qquad
\mbox{and}\\
w''(H)=\inf_{\underline{H}\subseteq H\subseteq\overline{H}}&\dfrac{\Big(w(\underline{H})+\varepsilon/2\Big)+\Big(w(\overline{H})+\varepsilon/2\Big)}{2}
\leq w(H)+\dfrac{\varepsilon}{2}.
}
Utilizing these two inequalities in \eq{700}, we get \eq{369}. This completes the proof.
\end{proof}
This research leaves many tempting questions open for further exploration. One of the directions is to study stability results for specific classes of graphs, such as trees, bipartite graphs, regular graphs, or complete graphs, where one may obtain tighter, case-specific norm bounds. Another promising direction is the generalization of the properties discussed here, such as approximate monotonicity, subadditivity, and convexity of graphs. Instead of a single scalar error term $\varepsilon>0$, one could associate a multidimensional function to model the error. This added flexibility may significantly broaden the scope of the theory and enhance its potential for real-world applications.

 \section*{Statements and Declarations}

\noindent\textbf{Funding.} 
The first author received financial support from the \textit{``Ipar a Veszprémi Mérnökképzésért'' Foundation} for the research presented in this work.

\vspace{6pt}
\noindent\textbf{Competing Interests.} 
The authors declare that there are no financial or non-financial competing interests relevant to the contents of this article.

\vspace{6pt}
\noindent\textbf{Ethics Approval.} 
Not applicable. This study does not involve human participants or animals.

\vspace{6pt}
\noindent\textbf{Consent to Participate.} 
Not applicable.

\vspace{6pt}
\noindent\textbf{Consent for Publication.} 
Not applicable.

\vspace{6pt}
\noindent\textbf{Data, Materials and/or Code Availability.} 
No datasets or code were generated or analysed during the current study.

\vspace{6pt}
\noindent\textbf{Authors’ Contributions.} 
The authors solely conceived the research idea, developed the theoretical framework, proved the results, and prepared the manuscript. They contributed equally in all aspects.


\end{document}